\documentclass[11pt]{article}
\usepackage{enumerate}
\usepackage{amssymb,a4wide,latexsym,makeidx,epsfig,fleqn}
\usepackage{amsthm}
\usepackage{amsmath}
\usepackage{enumerate}
\newtheorem{theorem}{Theorem}[section]

\newtheorem{lemma}[theorem]{Lemma}

\begin{document}
\textwidth 150mm \textheight 225mm
\title{Distance signless Laplacian spectral radius and Hamiltonian properties of graphs
\thanks{ Supported by
the National Natural Science Foundation of China (No.11171273)}}
\author{{Qiannan Zhou, Ligong Wang\footnote {Corresponding author.}}\\
{\small Department of Applied Mathematics, School of Science, Northwestern
Polytechnical University,}\\ 
{\small  Xi'an, Shaanxi 710072, People's Republic of China.}\\
{\small E-mail: qnzhoumath@163.com, lgwangmath@163.com}\\
 }
\date{}
\maketitle
\begin{center}
\begin{minipage}{120mm}
\vskip 0.3cm
\begin{center}
{\small {\bf Abstract}}
\end{center}
{\small In this paper, first, we establish a sufficient condition for a bipartite graph to be
Hamilton-connected. Furthermore, we also give two sufficient conditions on distance signless
Laplacian spectral radius for a graph to be Hamilton-connected and traceable from every vertex,
respectively. Last, we obtain a sufficient condition for a graph to be Hamiltonian in terms of
the distance signless Laplacian spectral radius of $G^{C}$.

 }
\end{minipage}
\end{center}

\section{Introduction }
\label{sec:ch6-introduction}

In this paper, we only consider finite simple undirected graphs. We denote by $m$ the edge
number of $G$, $\delta$ the minimum degree of $G$, $d_{G}(v)$ or simply $d(v)$ the degree of $v$
in $G$. We use $G[X,Y]$ to denote a bipartite graph
with bipartition $(X,Y)$. Let $K_{n}$ be a complete graph of order $n$ and $K_{m,n}$ be a
complete bipartite graph with two parts having $m,n$ vertices. For two disjoint
graphs $G_{1}$ and $G_{2}$, the union of $G_{1}$ and $G_{2}$, denoted by $G_{1}+G_{2}$, is
defined as $V(G_{1}+G_{2})=V(G_{1})\cup V(G_{2})$, $E(G_{1}+G_{2})=E(G_{1})\cup E(G_{2})$; and the join of $G_{1}$
and $G_{2}$, denoted by $G_{1}\vee G_{2}$, is defined as $V(G_{1}\vee G_{2})=V(G_{1})\cup V(G_{2})$, and
$E(G_{1}\vee G_{2})=E(G_{1}\vee G_{2})\cup \{xy: x\in V(G_{1}),y\in V(G_{2})\}$. The union of $k$ disjoint
copies of the same graph $G$ is denoted by $kG$. Write $K_{n-1}+e$ for the complete graph on
$n-1$ vertices with a pendant edge, and $K_{n-1}+v$ for the complete graph on $n-1$ vertices
together with an isolated vertex. The complement $G^{C}$ of $G$ is the graph on $V(G)$ with edge
set $[V]^{2}\setminus E(G)$.

The distance between $u$ and $v$ in $G$, denoted by $d_{G}(u,v)$, is the length of a shortest
path from $u$ to $v$. The transmission $Tr(u)$ of a vertex $u$ is defined to be the sum of
distances from $u$ to all other vertices in $G$, i.e., $Tr(u)=\sum_{v\in V(G)}d_{G}(u,v)$. A
graph $G$ is said to be transmission regular if $Tr(u)$ is a constant for each $u\in V(G)$. The
transmission of $G$ is the sum of distances between every pair of vertices of $G$. We denote it
by $\sigma(G)$. Obviously, we have $\sigma(G)=\frac{1}{2}\sum_{u\in V(G)}Tr(u)$. The distance
matrix of $G$, denoted by $\mathcal{D}(G)$, is a symmetric real matrix with $(i,j)$-entry being
$d(v_{i},v_{j})$. It is easy to
see that $Tr(v_{i})$ is the sum of $i$-th row of $\mathcal{D}(G)$. Let $Tr(G)$ be the diagonal
matrix of the vertex transmissions in $G$. Distance signless Laplacian matrix of graph $G$ is
defined as $Q_{D}(G)=Tr(G)+\mathcal{D}(G)$. The largest eigenvalue of $Q_{D}(G)$, denoted by
$\rho_{D}(G)$, is called to be the distance signless Laplacian spectral radius of $G$.

Hamilton cycle (path) is a cycle (path) that passes through all the vertices of a graph. A
graph is Hamiltonian (traceable) if it contains a Hamilton cycle (Hamilton path). And a graph is
Hamilton-connected if every two vertices of $G$ are connected by a Hamilton path. A graph is
traceable from a vertex $x$ if it has a Hamilton $x$-path.

Determining whether a given graph is Hamiltonian or not is an old problem in graph theory. This
problem is proved to be an NP-hard problem \cite{RMKa}. Many graph theorists are interested in
finding sufficient conditions for Hamilton cycles in graphs for a long time. In recent years,
graph theorists tried to use spectral graph theory to solve this problem. Of course, there are
many sufficient conditions on spectral radius or the signless Laplacian spectral radius for a
graph to be Hamiltonian, traceable or Hamilton-connected. In 2003, Krivelevich and Sudakov first
proposed a sufficient condition on the spectrum of the adjacency matrix for a regular graph to be
Hamiltonian, where the graphs satisfying the given condition are pseudo-random. Some other
spectral conditions for Hamilton cycles and paths in graphs have been given in
\cite{SBFC,VCh,Niki,RJG1,RJG2,MKBS,Ning,RLWC,YuFa,ZhouB}. In this paper, we mainly consider the
relationship between distance signless Laplacian spectral radius and the Hamiltonian properties
of graphs. In other words, we try to use distance signless Laplacian spectral radius to judge
whether a graph is Hamilton-connected or not. And we also give three sufficient conditions in
terms of distance signless Laplacian spectral radius of $G$ and one sufficient condition in
terms of distance signless Laplacian spectral radius of $G^{C}$.

\section{Lemmas and Results}
\label{sec:ch-sufficient}

Let $H_{t,n-t}$ $(t\geq 1)$ be a bipartite graph obtained from $K_{n,n-t}$ by adding $t$
vertices which are adjacent to $t$ common vertices with degree $n-t$ in $K_{n,n-t}$,
respectively. As for Lemma \ref{le:c1}, Moon and Moser in \cite{JMLM} obtained the strict
inequality. Ferrara, Jacobson and Powell in \cite{FMJM} characterized maximal nonhamiltonian
bipartite graphs.

\noindent\begin{lemma}\label{le:c1}(\cite{JMLM,FMJM}). Let $G=G[X,Y]$ be a bipartite graph with minimum degree $\delta\geq t$ $(t\geq 1)$ and $m$ edges, where $|X|=|Y|=n\geq 2t$. If
$$m\geq n^{2}-tn+t^{2},$$
then $G$ is Hamiltonian unless $G=H_{t,n-t}$.
\end{lemma}

\noindent\begin{lemma}\label{le:c2}(\cite{Rzhou}). Let $G=G[X,Y]$ be a bipartite graph on
$n$ vertices, then
$$\rho_{D}(G)\geq \left\{\begin{array}{ll}
  3n-4, & \mbox{ if n is even }\\
  \frac{5n-8+\sqrt{n^{2}+8}}{2}, & \mbox{ if n is odd }\\
\end{array}
\right.$$
\end{lemma}

\noindent\begin{theorem}\label{th:c3} Let $G[X,Y]$ be a bipartite graph with minimum degree $\delta\geq t$ $(t\geq 1)$ and $m$ edges, where $|X|=|Y|=n\geq 2t$. If
$$\rho_{D}(G)\leq m-n^{2}+(t+6)n-(t^{2}+4),$$
then $G$ is Hamiltonian unless $G=H_{t,n-t}$.

\end{theorem}

\noindent {\bf Proof.} Because $|X|=|Y|=n\geq 2t$, $G$ is a bipartite graph with order $2n$.
By Lemma \ref{le:c2}, we have $\rho_{D}(G)\geq 6n-4$. Combining with the conditions of Theorem
\ref{th:c3}, we get
$$6n-4\leq \rho_{D}(G) \leq m-n^{2}+(t+6)n-(t^{2}+4),$$
then $m\geq n^{2}-tn+t^{2}$. By Lemma \ref{le:c1}, we obtain that $G$ is
Hamiltonian unless $G=H_{t,n-t}$.

The proof is complete.
\hfill$\blacksquare$

\noindent\begin{lemma}\label{le:c4}(\cite{Rzhou}). Let $G$ be a connected graph on $n$ vertices,
then $$\rho_{D}(G)\geq \frac{4\sigma(G)}{n},$$ with equality holds if and only if $G$ is transmission
regular.
\end{lemma}

Let $\mathbb{NP}_{1}=\{K_{3}\vee (K_{n-5}+2K_{1}), K_{6}\vee 6K_{1}, K_{4}\vee (K_{2}+3K_{1}),
5K_{1}\vee K_{5}, K_{4}\vee (K_{1,4}+K_{1}), K_{4}\vee (K_{1,3}+K_{2}),
K_{3}\vee K_{2,5}, K_{4}\vee 4K_{1}, K_{3}\vee (K_{1}+K_{1,3}), K_{3}\vee (K_{1,2}+K_{2}),
K_{2}\vee K_{2,4}\}$.

\noindent\begin{lemma}\label{le:c5}(\cite{zhou2}). Let $G$ be a connected graph on $n\geq 5$
vertices and $m$ edges with minimum degree $\delta\geq 3$. If $$m\geq \dbinom{n-2}{2}+6,$$ then
$G$ is Hamilton-connected unless $G\in \mathbb{NP}_{1}$.
\end{lemma}

$G$ is a graph of order $n$, if $X$ is an eigenvector of $Q_{D}(G)$ corresponding to eigenvalue
$\rho$, then there is a 1-1 map $\varphi$ from $V(G)$ to the entries of $X$, simply written
as $X_{u}=\varphi(u)$ for each $u\in V(G)$. $X_{u}$ is also called the value of $u$ given by $X$.
We can find that
\begin{equation}\label{eq:c1}
[\rho-Tr(u)]X_{u}=\sum_{v\in V(G)}d_{G}(u,v)X_{v},
\end{equation}
for each $u\in V(G)$.

\noindent\begin{theorem}\label{th:c6} Let $G$ be a connected graph on $n\geq 5$ vertices with
minimum degree $\delta \geq 3$. If $$\rho_{D}(G)\leq \frac{2n^{2}+6n-36}{n},$$ then $G$ is
Hamilton-connected.
\end{theorem}

\noindent {\bf Proof.} Let $v\in V(G)$, $d(v)$ be the degree of $v$ in $G$. Then
$$Tr(v)\geq d(v)\cdot 1+(n-1-d(v))\cdot 2=2(n-1)-d(v),$$
with equality holds if and only if the maximum distance between $v$ and other vertices in $G$
is at most 2. So,
$$\sigma(G)=\frac{1}{2}\sum_{v\in V(G)}Tr(v)\geq \frac{1}{2}\sum_{v\in V(G)}[2(n-1)-d(v)]=n(n-1)-m,$$
with equality holds if and only if the maximum distance between $v$ and other vertices in $G$
is at most 2. By Lemma \ref{le:c4}, we get
$$\rho_{D}(G)\geq \frac{4\sigma(G)}{n}\geq 4(n-1)-\frac{4m}{n}.$$
Then by the conditions of Theorem \ref{th:c6}, we have
$$4(n-1)-\frac{4m}{n}\leq \rho_{D}(G)\leq \frac{2n^{2}+6n-36}{n},$$
then $m\geq \dbinom{n-2}{2}+6$. Suppose that $G$ is not Hamilton-connected, by Lemma \ref{le:c5},
we have $G\in \mathbb{NP}_{1}$. By direct calculation (see Table 1), we obtain that all graphs
in $\mathbb{NP}_{1}$ except $G=K_{3}\vee (K_{n-5}+2K_{1})$ satisfy
$\rho_{D}(G)\geq \frac{2n^{2}+6n-36}{n}$, so we can get a contradiction.

\begin{table} \caption {Some data of some graphs in $\mathbb{NP}_{1}$.}
\centering

\vskip1mm
\footnotesize
\begin{tabular}{p{3.5cm}p{2cm}p{1.5cm}p{1.5cm}}
%\B   begin{tabular}{cccc}
\hline
$G$&$\rho_{D}(G)$&$n$&$\frac{2n^{2}+6n-36}{n}$\tabularnewline
\hline
$K_{6}\vee 6K_{1}$            &28.8102   &12            &27\\
$K_{4}\vee (K_{2}+3K_{1})$    &21.2319   &9             &20\\
$5K_{1}\vee K_{5}$            &23.4031   &10            &22.4\\
$K_{4}\vee (K_{1,4}+K_{1})$   &23.8062   &10            &22.4\\
$K_{4}\vee (K_{1,3}+K_{2})$   &23.5751   &10            &22.4\\
$K_{3}\vee K_{2,5}$           &23.5751   &10            &22.4\\
$K_{4}\vee 4K_{1}$            &18        &8             &17.5\\
$K_{3}\vee (K_{1,3}+K_{1})$   &18.5208   &8             &17.5\\
$K_{3}\vee (K_{1,2}+K_{2})$   &18.2789   &8             &17.5\\
$K_{2}\vee K_{2,4}$           &18.2381   &8             &17.5\\
\hline
\end{tabular}
\end{table}

As for $G=K_{3}\vee (K_{n-5}+2K_{1})$, let $X=(x_{1},x_{2},\ldots,x_{n})^{T}$ be the eigenvector
corresponding to $\rho$. By (1), all vertices of transmission $n-1$ have the same values given by
$X$, say $X_{1}$; all vertices of transmission $2n-5$ have the same values given by $X$, say $X_{2}$. Denote
by $X_{3}$ the values of the vertices of transmission $n+1$ given by $X$. Assume
$\tilde{X}=(X_{1},X_{2},X_{3})^{T}$. Hence, by (1), we have
$$(\rho-(2+2+(n-5)))X_{1}=2X_{1}+2X_{2}+(n-5)X_{3},$$
$$(\rho-(3+2+2\cdot(n-5)))X_{2}=3X_{1}+2X_{2}+2\cdot(n-5)X_{3},$$
$$(\rho-(3+2\cdot2+(n-6)))X_{3}=3X_{1}+2\cdot2X_{2}+(n-6)X_{3}.$$
Transform the above equations into a matrix equation $(\rho I-A)\tilde{X}=0$, we get
\begin{displaymath}
A=\left(
  \begin{array}{cccccc}
      n+1&          2&           n-5&\\
        3&       2n-3&         2n-10&\\
        3&          4&          2n-5&\\
  \end{array}
\right).
\end{displaymath}
Thus, $\rho_{D}(G)$ is the largest root of the following equation
$$\rho^{3}-(5n-7)\rho^{2}+(8n^{2}-31n+56)\rho-4n^{3}+26n^{2}-82n+80=0.$$
Let $f(x)=x^{3}-(5n-7)x^{2}+(8n^{2}-31n+56)x-4n^{3}+26n^{2}-82n+80$, then
$f'(x)=3x^{2}-2(5n-7)x+8n^{2}-31n+56$. Let $f'(x)=0$, we get two roots $x_{1}$ and $x_{2}$, such that
$f'(x_{1})=f'(x_{2})=0$, where
$$x_{1}=\frac{5n-7-\sqrt{n^{2}+23n-119}}{3}, x_{2}=\frac{5n-7+\sqrt{n^{2}+23n-119}}{3}.$$
Consider $f(\frac{2n^{2}+6n-36}{n})=-\frac{8(n^{5}-6n^{4}-70n^{3}+954n^{2}-4050n+5832)}{n^{3}}< 0$
for $n\geq 5$ and $\frac{2n^{2}+6n-36}{n}> x_{2}$ for $m\geq 5$, which implies $\rho_{D}(G)> \frac{2n^{2}+6n-36}{n}$, a contradiction.

The proof is complete.
\hfill$\blacksquare$

Let $\mathbb{NP}_{2}=\{K_{2}\vee (K_{n-4}+2K_{1}), K_{5}\vee 6K_{1}, K_{3}\vee (K_{2}+3K_{1}),
5K_{1}\vee K_{4}, K_{3}\vee (K_{1,4}+K_{1}), K_{3}\vee (K_{1,3}+K_{2}),
K_{2}\vee K_{2,5}, K_{3}\vee 4K_{1}, K_{2}\vee (K_{1}+K_{1,3}), K_{2}\vee (K_{1,2}+K_{2}),
K_{1}\vee K_{2,4}\}$.

\noindent\begin{lemma}\label{le:c7}(\cite{zhou2}). Let $G$ be a connected graph on $n\geq 4$
vertices and $m$ edges with minimum degree $\delta\geq 2$. If $$m\geq \dbinom{n-2}{2}+4,$$ then
$G$ is traceable from every vertex unless $G\in \mathbb{NP}_{2}$.
\end{lemma}

\noindent\begin{theorem}\label{th:c8} Let $G$ be a connected graph on $n\geq 4$ vertices with
minimum degree $\delta \geq 2$. If $$\rho_{D}(G)\leq \frac{2n^{2}+6n-28}{n},$$ then $G$ is
traceable from every vertex.
\end{theorem}

\noindent {\bf Proof.} From the proof of Theorem \ref{th:c6}, we have
$\rho_{D}(G)\geq 4(n-1)-\frac{4m}{n}$. By the condition of Theorem \ref{th:c8}, we get
$$4(n-1)-\frac{4m}{n}\leq \rho_{D}(G)\leq \frac{2n^{2}+6n-28}{n},$$
then $m\geq \dbinom{n-2}{2}+4$. Suppose that $G$ is not traceable from every vertex, by
Lemma \ref{le:c7}, we obtain $G\in \mathbb{NP}_{2}$. By direct calculation (see Table 2), we
obtain that all graphs in $\mathbb{NP}_{2}$ except $G=K_{2}\vee (K_{n-4}+2K_{1})$ satisfy
 $\rho_{D}(G)> \frac{2n^{2}+6n-28}{n}$, so we can get a contradiction.

\begin{table} \caption {Some data of some graphs in $\mathbb{NP}_{2}$.}
\centering

\vskip1mm
\footnotesize
\begin{tabular}{p{3.5cm}p{2cm}p{1.5cm}p{1.5cm}}
%\begin{tabular}{cccc}
\hline
$G$&$\rho_{D}(G)$&$n$&$\frac{2n^{2}+6n-28}{n}$\tabularnewline
\hline
$K_{5}\vee 6K_{1}$            &27.2621   &11            &25.455\\
$K_{3}\vee (K_{2}+3K_{1})$    &19.6847   &8             &18.5\\
$5K_{1}\vee K_{4}$            &21.8443   &9             &20.889\\
$K_{3}\vee (K_{1,4}+K_{1})$   &22.0660   &9             &20.889\\
$K_{3}\vee (K_{1,3}+K_{2})$   &22.0083   &9             &20.889\\
$K_{2}\vee K_{2,5}$           &22.0120   &9             &20.889\\
$K_{3}\vee 4K_{1}$            &16.4244   &7             &16\\
$K_{2}\vee (K_{1,3}+K_{1})$   &16.9667   &7             &16\\
$K_{2}\vee (K_{1,2}+K_{2})$   &16.6974   &7             &16\\
$K_{1}\vee K_{2,4}$           &16.6569   &7             &16\\
\hline
\end{tabular}
\end{table}

For $G=K_{2}\vee (K_{n-4}+2K_{1})$, let $X=(x_{1},x_{2},\ldots,x_{n})^{T}$ be the eigenvector
corresponding to $\rho$. By (1), all vertices of transmission $n-1$ have the same values given by
$X$, say $X_{1}$; all vertices of transmission $2n-4$ have the same values given by $X$, say $X_{2}$. Denote
by $X_{3}$ the values of the vertices of transmission $n+1$ given by $X$. Assume
$\tilde{X}=(X_{1},X_{2},X_{3})^{T}$. Hence, by (1), we have
$$(\rho-(1+2+(n-4)))X_{1}=X_{1}+2X_{2}+(n-4)X_{3},$$
$$(\rho-(2+2+2\cdot(n-4)))X_{2}=2X_{1}+2X_{2}+2\cdot(n-4)X_{3},$$
$$(\rho-(2+2\cdot2+(n-5)))X_{3}=2X_{1}+2\cdot2X_{2}+(n-5)X_{3}.$$
Transform the above equations into a matrix equation $(\rho I - B)\tilde{X}=0$, we get

\begin{displaymath}
B=\left(
  \begin{array}{cccccc}
        n&          2&           n-4&\\
        2&       2n-2&          2n-8&\\
        2&          4&         2n-4&\\
  \end{array}
\right).
\end{displaymath}
Thus, $\rho_{D}(G)$ is the largest root of the following equation:
$$\rho^{3}-(5n-6)\rho^{2}+(8n^{2}-28n+44)\rho-4n^{3}+24n^{2}-68n+64=0.$$
Let $g(x)=x^{3}-(5n-6)x^{2}+(8n^{2}-28n+44)x-4n^{3}+24n^{2}-68n+64=0$, then
$g'(x)=3x^{2}-2(5n-6)x+8n^{2}-28n+44$. Let $g'(x)=0$, we have two values $x_{1}$ and $x_{2}$, such that
$g'(x_{1})=g'(x_{2})=0$, where
$$x_{1}=\frac{5n-6-\sqrt{n^{2}+24n-96}}{3}, x_{2}=\frac{5n-6+\sqrt{n^{2}+24n-96}}{3}.$$
Consider $g(\frac{2n^{2}+6n-28}{n})=-\frac{8(n^{5}-4n^{4}-67n^{3}+686n^{2}-2352n+2744)}{n^{3}}<0$
for $n\geq 4$ and $\frac{2n^{2}+6n-28}{n}> x_{2}$ for $n\geq 4$, which implies $\rho_{D}(G)> \frac{2n^{2}+6n-28}{n}$, we can get a contradiction.

The proof is complete.
\hfill$\blacksquare$

\noindent\begin{lemma}\label{le:c9}(\cite{Niki}). Let $G$ be a graph on $n$ vertices and $m$
edges. If
$$m\geq \dbinom{n-1}{2},$$
then $G$ contains a Hamilton path unless $G=K_{n-1}+v$. If the inequality is strict, then $G$
contains a Hamilton cycle unless $G=K_{n-1}+e$.
\end{lemma}

\noindent\begin{theorem}\label{th:c10} Let $G$ be a graph on $n$ vertices, $m$ edges and
$\rho_{D}(G^{C})$ be the distance signless Laplacian spectral radius of its complement. If
$$\rho_{D}(G^{C})\leq \frac{3n^{2}-n+10m-2}{2n},$$
then $G$ contains a Hamilton path unless $G=K_{n-1}+v$. If the inequality is strict, then
$G$ contains a Hamilton cycle unless $G=K_{n-1}+e$.
\end{theorem}

\noindent {\bf Proof.} Because $\sigma(G^{C})=\frac{1}{2}\sum_{v\in V(G)}Tr_{G^{C}}(v)
\geq \frac{1}{2}\sum_{v\in V(G^{C})}[(n-d(v))+2d(v)]
=\frac{1}{2}\sum_{v\in V(G^{C})}[n+d(v)]=\frac{1}{2}n(n-1)+m$. By Lemma \ref{le:c4}, we have
$$\rho_{D}(G^{C})\geq \frac{4\sigma(G^{C})}{n}\geq 2(n-1)+\frac{4m}{n}.$$
Combining with the condition of Theorem \ref{th:c10}, we get
$$2(n-1)+\frac{4m}{n}\leq \rho_{D}(G^{C}) \leq \frac{3n^{2}-n+10m-2}{2n},$$
so $m\geq \dbinom{n-1}{2}$. Then by Lemma \ref{le:c9}, we can obtain the conclusion.

The proof is complete.
\hfill$\blacksquare$

\end{document}